\documentclass[english]{article}
\usepackage[T1]{fontenc}
\usepackage[latin9]{inputenc}
\usepackage{mathrsfs}
\usepackage{amsmath}
\usepackage{amssymb}
\usepackage{stackrel}

\makeatletter
\newcommand{\lyxaddress}[1]{
\par {\raggedright #1
\vspace{1.4em}
\noindent\par}
}

\makeatother

\usepackage{babel}
\begin{document}

\title{Some results on multi vector space}

\author{Moumita Chiney$^{1}$, S. K. Samanta$^{2}$}

\maketitle

\lyxaddress{Department of Mathematics, Visva-Bharati, Santiniketan-731235. \\
Email: $^{1}$moumi.chiney@gmail.com, $^{2}$ syamal.samanta@visva-bharati.ac.in,
syamal\_123@yahoo.co.in}
\begin{abstract}
In the present paper, a notion of M-basis and multi dimension of a
multi vector space is introduced and some of its properties are studied.
\end{abstract}

\section{\noindent Introduction}

Theory of Multisets is an important generalization of classical set
theory which has emerged by violating a basic property of classical
sets that an element can belong to a set only once. Synonymous terms
of multisets viz. list, heap, bunch, bag, sample, weighted set, occurrence
set and fireset (finitely repeated element set) are used in different
contexts but conveying the same idea. It is a set where an element
can occur more than once. Many authors like Yager \cite{yager}, Miyamoto
\cite{miyamoto,miyamoto1}, Hickman \cite{hickman}, Blizard \cite{bli},
Girish and John \cite{girish,girish-1}, Chakraborty \cite{chakraborty}
etc. have studied the properties of multisets. Multisets are very
useful structures arising in many areas of mathematics and computer
science such as database queries, multicriteria decision making, knowledge
representation in data based systems, biological systems and membrane
computing etc. \cite{delago,delago2,klausner,kosters,li,miyamoto1,mumick,paun}.
Again the theory of vector space is one of the most important algebraic
structures in modern mathematics and this has been extended in different
setting such as fuzzy vector space \cite{katsaras,lubzonok,shi},
intuitionistic fuzzy vector space , soft vector space \cite{sujoy}
etc. In \cite{chiney}, we introduced a notion of multi vector space
and studied some of its basic properties. As a continuation of our
earlier paper \cite{chiney}, here we have attempted to formulate
the concept of basis and dimension of multi vector space and to study
their properties.

\section{\noindent Preliminaries}

In this section definition of a multiset (mset in short) and some
of its properties are given. Unless otherwise stated, $X$ will be
assumed to be an initial universal set and $MS(X)$ denotes the set
of all mset over $X.$ $\newline$ \\
$\mathbf{Definition\:2.1}$ \cite{girish-1} An mset $M$ drawn from
the set $X$ is represented by a count function $C_{M}:X\rightarrow N$
where $N$ represents the set of non negative integers.

Here $C_{M}(x)$ is the number of occurrence of the element $x$ in
the mset $M.$ The presentation of the mset $M$ drawn from $X=\{x_{1},x_{2},....,x_{n}\}$
will be as $M=\{x_{1}/m_{1},x_{2}/m_{2},....,x_{n}/m_{n}\}$ where
$m_{i}$ is the number of occurrences of the element $x_{i},$ $i=1,2,...,n$
in the mset $M.$

Also here for any positive integer $\omega,$ $[X]^{\omega}$ is the
set of all msets whose elements are in $X$ such that no element in
the mset occurs more than $\omega$ times and it will be referred
to as mset spaces.$\newline$ \\
The algebraic operations of msets are considered as in \cite{girish-1}.$\newline$\\
$\mathbf{Definition\:2.2}$ \cite{nazmul} Let $M$ be a mset over
a set $X$. Then a set $M_{n}=\{x\in X:C_{M}(x)\geq n\},$ where $n$
is a natural number, is called $n-$ level set of $M.$ $\newline$\\
$\mathbf{Proposition\:2.3}$ \cite{nazmul} Let $A,B$ be msets over
$X$ and $m,n\in\mathbb{N}$. \\
$(1)$ If $A\subseteq B,$ then $A_{n}\subseteq B_{n};$\\
$(2)$ If $m\leq n,$ then $A_{m}\supseteq A_{n};$\\
$(3)$ $(A\cap B)_{n}=A_{n}\cap B_{n};$\\
$(4)$ $(A\cup B)_{n}=A_{n}\cup B_{n};$\\
$(5)$ $A=B$ iff $A_{n}=B_{n},\forall n\in\mathbb{N}.$ $\newline$\\
$\mathbf{Definition\:2.4}$ \cite{nazmul} Let $P\subseteq X.$ Then
for each $n\in\mathbb{N},$ we define a mset $nP$ over $X,$ where
$C_{nP}(x)=n,\forall x\in P.$$\newline$\\
$\mathbf{Definition\:2.5}$ \cite{nazmul} Let $X$ and $Y$ be two
nonempty sets and $f:X\rightarrow Y$ be a mapping. Then \\
$(1)$ the image of a mset $M\in[X]^{\omega}$ under the mapping $f$
is denoted by $f(M)$ or $f[M]$, where \\
\[
C_{f(M)}(y)=\begin{cases}
\underset{f(x)=y}{\vee} & C_{M}(x)\: if\: f^{-1}(y)\neq\phi\\
0 & otherwise
\end{cases}
\]
$(2)$ the inverse image of a mset $N\in[Y]^{\omega}$ under the mapping
$f$ is denoted by $f^{-1}(N)$ or $f^{-1}[N]$ where $C_{f^{-1}(N)}(x)=C_{N}[f(x)]$.$\newline$\\
The properties of functions, which are used in this paper, are as
in \cite{nazmul}.$\newline$

\noindent $\mathbf{Definition\:2.6}$ \cite{chiney} Let $A_{1},A_{2},...,A_{n}$
$\in\left[X\right]^{\omega}.$ Then we define $A_{1}+A_{2}+...+A_{n}$
as follows:\\
$C_{A_{1}+A_{2}+...+A_{n}}(x)=\vee\{C_{A_{1}}(x_{1})\wedge C_{A_{2}}(x_{2})\wedge...\wedge C_{A_{n}}(x_{n}):x_{1},x_{2},...,x_{n}\in X\: and\; x_{1}+x_{2}+...+x_{n}=x\}.$
\\
Let $\lambda\in K$ and $B\in\left[X\right]^{\omega}.$ Then $\lambda B$
is defined as follows:\\
$C_{\lambda B}(y)=\vee\{C_{B}(x):\lambda x=y\}.$$\newline$\\
$\mathbf{Lemma\:2.7}$ \cite{chiney} Let $\lambda\in K$ and $B\in\left[X\right]^{\omega}.$
Then \\
$(a)$ For $\lambda\neq0,$ $C_{\lambda B}(y)=C_{B}(\lambda^{-1}y),\forall y\in X.$
\\
For $\lambda=0,$ \\
$C_{\lambda B}(y)=\begin{cases}
0, & y\neq0,\\
\underset{x\in X}{sup}\: C_{B}(x), & y=0
\end{cases}$.\\
$(b)$ For all scalars $\lambda\in K$ and for all $x\in X,$ we have
$C_{\lambda B}(\lambda x)\geq C_{B}(x).$$\newline$

\noindent $\mathbf{Definition\:2.8}$ \cite{chiney} A multiset $V$
in $[X]^{\omega}$ is said to be a multi vector space or multi linear
space(in short mvector space) over the linear space $X$ if \\
$(i)$ $V+V\subseteq V;$\\
$(ii)$ $\lambda V\subseteq V,$ for every scalar $\lambda.$\\
We denote the set of all multi vector space over $X$ by $MV(X).$
$\newline$\\
$\mathbf{Remark\:2.9}$ \cite{chiney} For a multi vector space $V$
in $[X]^{\omega}$ , $V+V+.....n\: times=V,$ i. e., $nV=V.$ $\newline$\\
$\mathbf{Remark\:2.10}$ \cite{chiney} If $V\in MV(X)$ with $dim\: X=m,$
then $\mid C_{V}(X)\mid\leq m+1,$ where $\mid C_{V}(X)\mid$ represents
the cardinality of $C_{V}(X).$$\newline$\\
$\mathbf{Proposition\:2.11}$ \cite{chiney} (\textbf{Representation
theorem}) Let $V\in MV(X)$ with $dim\: X=m$ and range of $C_{V}=\{n_{0},n_{1},....,n_{k}\}\subseteq\{0,1,2,...,\omega\},$
$k\leq m,$ $n_{0}=C_{V}(\theta)$ and $\omega\geq n_{0}>n_{1}>...>n_{k}\geq0.$
Then there exists a nested collection of subspaces of $X$ as\\
$\{\theta\}\subseteq V_{n_{0}}\subsetneqq V_{n_{1}}\subsetneqq V_{n_{2}}\subsetneqq....\subsetneqq V_{n_{k}}=X$
such that $V=n_{0}V_{n_{0}}\cup n_{1}V_{n_{1}}\cup.....\cup n_{k}V_{n_{k}}.$
Also\\
$(1)$ If $n,m\in(n_{i+1},n_{i}],$ then $V_{n}=V_{m}=V_{n_{i}}.$
\\
$(2)$ If $n\in(n_{i+1},n_{i}]$ and $m\in(n_{i},n_{i-1}]$, then
$V_{n}\supsetneqq V_{m}.$ $\newline$\\
$\mathbf{Definition\:2.12}$ \cite{chiney} Let $X$ be a finite dimensional
vector space with $dim\: X=m$ and $V\in MV(X).$ Consider $Proposition\:2.11.$
Let $B_{n_{i}}$ be a basis on $V_{n_{i}},\: i=0,1,...,k$ such that\\
$B_{n_{0}}\subsetneqq B_{n_{1}}\subsetneqq B_{n_{2}}\subsetneqq...\subsetneqq B_{n_{k}}$.....$(iii)$\\
 Define a multi subset $\beta$ of $X$ by\\
$C_{\beta}(x)=\begin{cases}
\vee\{n_{i}:x\in B_{n_{i}}\}\\
0,\: otherwise
\end{cases}$\\
Then $\beta$ is called a multi basis of $V$ corresponding to $(iii).$
We denote the set of all multi bases of $V$ by $\mathcal{B}_{M}(V).$$\newline$\\
$\mathbf{Corollary\:2.13}$ \cite{chiney} Let $\beta$ be a multi
basis of $V$ obtained by $(iii).$ Then\\
$(1)$ If $n,m\in(n_{i+1},n_{i}]$, then $\beta_{n}=\beta_{m}=B_{n_{i}}.$\\
$(2)$ If $n\in(n_{i+1},n_{i}]$ and $m\in(n_{i},n_{i-1}]$, then
$\beta_{n}\supsetneqq\beta_{m}.$\\
$(3)$ $\beta_{n}$ is a basis of $V_{n}$, for all $n\in\{1,2,....,\omega\}$.
$\newline$

\section{Some results on multi vector space}

$\mathbf{Lemma\:3.1}$ Let $s,t\in\mathbb{R}$ and $A,\: A_{1}$ and
$A_{2}$ be multisets on a vector space $X.$ Then \\
$(1)$ $s.(t.A)=t.(s.A)=(st).A$ and\\
$(2)$ $A_{1}\leq A_{2}\Rightarrow t.A_{1}\leq t.A_{2}.$$\newline$\\
$\mathbf{Proposition\:3.2}$ Let $V\in MV(X).$ Then $x\in X,\; a\neq0\Rightarrow C_{V}(ax)=C_{V}(x).$$\newline$\\
$\mathbf{Proposition\:3.3}$ Let $V\in MV(X)$ and $u,v\in X$ such
that $C_{V}(u)>C_{V}(v)$. Then $C_{V}(u+v)=C_{V}(v).$ $\newline$\\
$\mathbf{Proposition\:3.4}$ Let $V\in MV(X)$ and $v,w\in X$ with
$C_{V}(v)\neq C_{V}(w)$. Then $C_{V}(v+w)=C_{V}(v)\wedge C_{V}(w).$
$\newline$\\

\section{Multi linear independence}

$\mathbf{Definition\:4.1}$ Let $V\in MV(X)$ and $dim\: X=n.$ We
say that a finite set of vectors $\{x_{i}\}_{i=1}^{n}$ is multi linearly
independent in $V$ if and only if $\{x_{i}\}_{i=1}^{n}$ is linearly
independent in $X$ and for all $\{a_{i}\}_{i=1}^{n}\subset\mathbb{R}$
with $a_{i}\neq0$ , $C_{V}(\stackrel[i=1]{n}{\sum}a_{i}x_{i})$\\
$=\stackrel[i=1]{n}{\wedge}C_{V}(a_{i}x_{i}).$ $\newline$ \\
The following example shows that every linearly independent set is
not multi linearly independent.\\
$\mathbf{Example\;4.2}$ Let $X=\mathbb{R}^{2}$ and $\omega=4.$
We define a multi vector space $C_{V}:X\rightarrow N$ by 

$C_{V}(x)=\begin{cases}
4,\: if\: x=(0,0)\\
2,\; if\; x=(0,a),a\neq0\\
1,\: otherwise.
\end{cases}$\\
If we take the vectors $x=(1,0)$ and $y=(-1,1),$ then they are linearly
independent but not multi linearly independent. As here $C_{V}(x)=C_{V}(y)=1$,
but $C_{V}(x+y)=2>(C_{V}(x)\wedge C_{V}(y))=1.$ $\newline$\\
$\mathbf{Proposition\:4.3}$ Let $V\in MV(X)$ and $dim\: X=m.$ Then
any set of vectors $\{x_{i}\}_{i=1}^{N}(N\leq m),$ which have distinct
counts is linearly and multi linearly independent.$\newline$ \\
$Proof.$ The proof follows by method of induction.$\newline$\\
$\mathbf{Note\;4.4}$ Converse of the above proposition is not true.
Let $X=\mathbb{R}^{2}$ and $\omega=6.$ We define a multi vector
space $C_{V}:X\rightarrow N$ by 

$C_{V}(x)=\begin{cases}
6,\: if\: x=(0,0)\\
1,\: otherwise.
\end{cases}$\\
Then we have $\{\theta\}=V_{6}\subset V_{1}=\mathbb{R}^{2}.$ Let
$e_{1}=(1,0),$ $e_{2}=(0,1).$ Then $\{e_{1},e_{2}\}$ are multi
linearly independent in $V$ , although, $C_{V}(e_{1})=C_{V}(e_{2})$.
$\newline$

\section{M-basis}

$\mathbf{Definition\:5.1}$ A M-basis for a multi vector space $V\in MV(X)$
is a basis of $X$ which is multi linearly independent in $V$. \\
We denote the set of all M-bases of $V$ by $\text{\ensuremath{\mathscr{B}}}(V).$$\newline$\\
$\mathbf{Lemma\:5.2}$ If $V\in MV(X)$ and $Y$ is a proper subspace
of $X,$ then for any $t\in X\setminus Y$ with $C_{V}(t)=sup[C_{V}(X\setminus Y)]$,
$C_{V}(t+y)=C_{V}(t)\wedge C_{V}(y),$ for all $y\in Y.$ $\newline$\\
$Proof.$ Since $\omega$ is finite, such a $t$ exists. Let $y\in Y.$
If $C_{V}(y)\neq C_{V}(t)$ then by $Proposition\:3.4,$ $C_{V}(t+y)=C_{V}(t)\wedge C_{V}(y).$
If $C_{V}(y)=C_{V}(t)$ then by $Definition\:2.8,$ $C_{V}(t+y)\geq C_{V}(t)\wedge C_{V}(y).$
Since $t+y\in X\setminus Y$ and $C_{V}(t)=sup[C_{V}(X\setminus Y)],$
we must have $C_{V}(t+y)\leq C_{V}(t)=C_{V}(y)$ and thus $C_{V}(t+y)=C_{V}(t)\wedge C_{V}(y).$
$\newline$\\
$\mathbf{Lemma\:5.3}$ Let $V\in MV(X)$, $Y$ be a proper subspace
of $X$ and $C_{V}\mid_{Y}=C_{V^{\prime}}.$ If $B^{*}$ is a M-basis
for $V^{\prime},$ then there exists $t\in X\setminus Y$ such that
$B^{+}=B^{*}\cup\{t\}$ is a M-basis for $W,$ where $C_{W}=C_{V}\mid_{\prec B^{+}\text{\ensuremath{\succ}}}$
and $\text{\ensuremath{\prec}}B^{+}\text{\ensuremath{\succ}}$ is
the vector space spanned by $B^{+}.$ $\newline$\\
$Proof.$ Pick $t\in X\setminus Y$ such that $C_{V}(t)=sup[C_{V}(X\setminus Y)].$
Then by $Lemma\;5.2,$ $B^{+}=B^{*}\cup\{t\}$ is a multi linearly
independent and hence a M-basis for $W,$ where $C_{W}=C_{V}\mid_{\prec B^{+}\text{\ensuremath{\succ}}}.$
$\newline$\\
$\mathbf{Proposition\:5.4}$ All multi vector spaces $V\in MV(X)$
with $dim\: X=m$ have M-basis.$\newline$\\
$Proof.$ The proof follows by mathematical induction. $\newline$\\
$\mathbf{Proposition\:5.5}$ Let $V\in MV(X)$ where $dim\: X=m$
and $C_{V}(X\setminus\{\theta\})=\{n_{0},n_{1},n_{2},...,n_{k}\},\; k\leq m$.
Then a basis $B=\{e_{1},e_{2},...,e_{m}\}$ of $X$  is a M-basis
for $V$ if and only if $B\cap V_{n_{i}}$ is a basis of $V_{n_{i}}$
for any $i=0,1,...,k.$ $\newline$\\
$\mathbf{Proof.}$ Let $\omega\geq n_{0}>n_{1}>....>n_{k}\geq0.$
Then $\{\theta\}\subsetneqq V_{n_{0}}\subsetneqq V_{n_{1}}\subsetneqq V_{n_{2}}\subsetneqq....\subsetneqq V_{n_{k}}=X$
. Let $B_{n_{i}}=B\cap V_{n_{i}},\; i=0,1,...,k.$ \\
First suppose that $B\cap V_{n_{i}}=B_{n_{i}}$ is a basis of $V_{n_{i}}$
for any $i=0,1,...,k.$ Then $B_{n_{0}}\subsetneqq B_{n_{1}}\subsetneqq......\subsetneqq B_{n_{k}}=B.$
Let $B_{n_{0}}=\{e_{n_{0}}^{1},e_{n_{0}}^{2},...,e_{n_{0}}^{j}\},j\leq m.$
Then $C_{V}(\stackrel[i=1]{j}{\sum}a_{i}e_{n_{0}}^{i})\geq\stackrel[i=1]{j}{\wedge}C_{V}(e_{n_{0}}^{i})=n_{0}$.
Since $n_{0}$ is the highest count,\\
 $C_{V}(\stackrel[i=1]{j}{\sum}a_{i}e_{n_{0}}^{i})=n_{0}=\stackrel[i=1]{j}{\wedge}C_{V}(e_{n_{0}}^{i}).$
Hence $B_{n_{0}}$ is multi linearly independent. \\
Next let $B_{n_{1}}=B_{n_{0}}\cup\{e_{n_{1}}^{1},e_{n_{1}}^{2},...,e_{n_{1}}^{s}\},j+s\leq m.$
Consider the sum\\
 $\stackrel[i=1]{j}{\sum}b_{i}e_{n_{0}}^{i}+\stackrel[i=1]{s}{\sum}c_{i}e_{n_{1}}^{i}$,
where some $c_{i}\neq0.$ Then $C_{V}(\stackrel[i=1]{j}{\sum}b_{i}e_{n_{0}}^{i}+\stackrel[i=1]{s}{\sum}c_{i}e_{n_{1}}^{i})$\\
$\geq\left(\underset{i\in J_{1}}{\wedge}C_{V}(e_{n_{0}}^{i})\right)\wedge\left(\underset{i\in J_{2}}{\wedge}C_{V}(e_{n_{1}}^{i})\right)$,
{[}where $J_{1}=\{i:b_{i}\neq0\},\; J_{2}=\{i:c_{i}\neq0\}${]} $=n_{1}.$\\
 If $C_{V}(\stackrel[i=1]{j}{\sum}b_{i}e_{n_{0}}^{i}+\stackrel[i=1]{s}{\sum}c_{i}e_{n_{1}}^{i})>n_{1},$
then \\
$C_{V}(\stackrel[i=1]{j}{\sum}b_{i}e_{n_{0}}^{i}+\stackrel[i=1]{s}{\sum}c_{i}e_{n_{1}}^{i})=n_{0}\Rightarrow\stackrel[i=1]{j}{\sum}b_{i}e_{n_{0}}^{i}+\stackrel[i=1]{s}{\sum}c_{i}e_{n_{1}}^{i}\in V_{n_{0}}\Rightarrow c_{i}=0,$
for all $i=1,2,..,s,$ a contradiction. Thus $B_{n_{1}}$ is multi
linearly independent. Proceeding in the similar way it can be proved
that $B_{n_{k}}=B$ is multi linearly independent and hence a M-basis
for $V.$\\
Conversely, let $B$ be a M-basis for $V.$ Then either $B_{n_{i}}=\phi$
or $B_{n_{i}}\neq\phi.$ \\
Let $B_{n_{i}}=\phi$ and $x\in V_{n_{i}}.$ Then obviously $B_{n_{j}}=\phi,j<i.$
Since $B$ is a basis of $X,$ there exists some $B^{\prime}\subseteq B$
such that $x=\underset{e_{j}\in B}{\sum}b_{j}e_{j},$ $b_{j}\neq0.$
Then $C_{V}(x)=\underset{e_{j}\in B}{\wedge}C_{V}(e_{j})\leq n_{i+1},$
a contradiction. So, $B_{n_{i}}\neq\phi.$\\
Then $B_{n_{0}}\subsetneqq B_{n_{1}}\subsetneqq......\subsetneqq B_{n_{k}}=B.$
Let $x\in V_{n_{i}}$ and $B_{n_{i}}$ is not a basis of $V_{n_{i}}$.
Choose $x=\underset{e_{i}\in B_{n_{i}}}{\sum}a_{i}e_{i}+\underset{e_{i}^{\prime}\notin B_{n_{i}}}{\sum}b_{i}e_{i}^{\prime}$,
for all $b_{i}\neq0$. \\
Now, $C_{V}(x)=C_{V}(\underset{e_{i}\in B_{n_{i}}}{\sum}a_{i}e_{i}+\underset{e_{i}^{\prime}\notin B_{n_{i}}}{\sum}b_{i}e_{i}^{\prime})$
\\
$=\left(\underset{e_{i}\in B_{n_{i}}^{\prime}}{\wedge}C_{V}(e_{i})\right)\wedge\left(\underset{e_{i}^{\prime}\notin B_{n_{i}}}{\wedge}C_{V}(e_{i}^{\prime})\right)$,
{[}where $B_{n_{i}}^{\prime}=\{e_{i}\in B_{n_{i}}:a_{i}\neq0\}${]}
$=\left(\underset{e_{i}^{\prime}\notin B_{n_{i}}}{\wedge}C_{V}(e_{i}^{\prime})\right)<n_{i},$
a contradiction to the fact that $x\in V_{n_{i}}.$ Thus $x=\underset{e_{i}\in B_{n_{i}}}{\sum}a_{i}e_{i}$
and $B_{n_{i}}$ is a basis of $V_{n_{i}}.$ \\
Hence proved. $\newline$\\
$\mathbf{Proposition\:5.6}$ Let $V$ be a multi vector space over
$X$ where $dim\: X=m.$ Then there is an one-to-one correspondence
between $\mathcal{B}_{M}(V)$ and $\text{\ensuremath{\mathscr{B}}}(V).$
$\newline$\\
$\mathbf{Proposition\;5.7}$ Let $V\in MV(X)$ with $dim\: X=m$ and
range of $C_{V}(X\setminus\{\theta\})=\{n_{0},n_{1},....,n_{k}\}\subseteq\{0,1,2,...,\omega\},$
$k\leq m.$ If a basis $B=\{e_{1},e_{2},...,e_{m}\}$ of $X$ is a
M-basis, then $C_{V}(B)=\{n_{0},n_{1},....,n_{k}\}.$ $\newline$\\
$\mathbf{Remark\;5.8}$ Converse of the above $Proposition$ is not
true. For example, suppose $X=\mathbb{R}^{4},\;\omega=5.$ Define
multi vector space $V$ with count functions $C_{V}$ as follows:\\
$C_{V}((0,0,0,0))=5$; $C_{V}((0,0,0,\mathbb{R}\setminus\{0\}))=5;$
$C_{V}((0,0,\mathbb{R}\setminus\{0\},\mathbb{R}))=5$, $C_{V}((0,\mathbb{R}\setminus\{0\},\mathbb{R},\mathbb{R}))=2,C_{V}(\mathbb{R}^{4}\setminus(0,\mathbb{R},\mathbb{R},\mathbb{R}))=2.$
\\
Then $B=\{(0,0,0,1),(-1,1,1,1),(1,-1,1,1),(1,1,-1,1)\}$ is a basis
of $\mathbb{R}^{4}$ and $C_{V}(B)=\{2,5\}=C_{V}(\mathbb{R}^{4})$.
But $B$ is not a M-basis as $B$ is not multi linearly independent.
In fact, $C_{V}((-1,1,1,1))=C_{V}(1,-1,1,1))=2.$ But $C_{V}((-1,1,1,1)+(1,-1,1,1))=C_{V}((0,0,1,1))=5>[C_{V}((-1,1,1,1))\wedge C_{V}((-1,1,1,1))]=2.$$\newline$\\
$\mathbf{Definition\;5.9}$ Let $V\in MV(X)$ with $dim\: X=m$, range
of $C_{V}(X\setminus\{\theta\})=\{n_{0},n_{1},$$...,n_{k}\}$$\subseteq\{0,1,2,...,\omega\},$
$k\leq m$ and $B_{0}$ be any M-basis for $V.$ Then $C_{V}(B_{0})=\{n_{0},n_{1},....,n_{k}\}.$
We define multi index of a multi M-basis $B_{0}$ with respect to
$V$ by $[B_{0}]_{M}=\{r_{i}:$ $r_{i}$ is the number of element
of $B_{0}$ taking the value $n_{i}$$\}.$$\newline$\\
$\mathbf{Proposition\;5.10}$ For a multi vector space $V$, multi
index of M-basis with respect to $V$ is independent of M-basis. $\newline$\\
$\mathbf{Proof.}$ Let $V\in MV(X)$ with $dim\: X=m$, range of $C_{V}(X\setminus\{\theta\})=\{n_{0},n_{1},....,n_{k}\}$\\
$\subseteq\{0,1,2,...,\omega\},$ $k\leq m$ and $\omega\geq n_{0}>n_{1}>...>n_{k}\geq0.$
Then for any two M-bases $B_{0},B_{0}^{\prime}$ of $V,$ $C_{V}(B_{0})=C_{V}(B_{0}^{\prime})=\{n_{0},n_{1},....,n_{k}\}.$
Let $[B_{0}]_{M}=\{r_{i}\}$ and $[B_{0}^{\prime}]_{M}=\{r_{i}^{\prime}\}$.
Now, $\mid B_{0}\cap V_{n_{i}}\mid=\stackrel[j=0]{i}{\sum}r_{j}$
and $\mid B_{0}^{\prime}\cap V_{n_{i}}\mid=\stackrel[j=0]{i}{\sum}r_{j}^{\prime}$,
for $i=0,1,2,...,k.$ As $B_{0}\cap V_{n_{i}}$ and $B_{0}^{\prime}\cap V_{n_{i}}$
are both basis of $V_{n_{i}}$, $\mid B_{0}\cap V_{n_{i}}\mid$ \\
$=\mid B_{0}^{\prime}\cap V_{n_{i}}\mid,$ for all $i=0,1,2,...,k.$
Thus $[B_{0}]_{M}=[B_{0}^{\prime}]_{M}$.$\newline$\\
$\mathbf{Note\;5.11}$ As multi index of M-basis with respect to a
multi vector space $V$ is independent of M-basis, we can use only
the term multi index of a multi vector space $V.$ $\newline$\\
$\mathbf{Definition\;5.12}$ Let $V\in MV(X)$ with $dim\: X=m$,
$C_{V}(X)=\{n_{0},n_{1},....,n_{k}\}$\\
$\subseteq\{0,1,2,...,\omega\},$ $k\leq m$ and $B$ be any basis
for $X.$ We define index of a basis $B$ with respect to $V$ by
$[B]=\{r_{i}:\; n_{i}$ $r_{i}$ is the number of element of $B$
taking the value $n_{i}$$\}.$$\newline$\\
$\mathbf{Proposition\;5.13}$ Let $V\in MV(X)$ with $dim\: X=m$,
$C_{V}(X\setminus\{\theta\})=\{n_{0},n_{1},....$\\
$,n_{k}\}$$\subseteq\{0,1,2,...,\omega\},$ $k\leq m$ and $B$ be
any basis of $X$ with $C_{V}(B)=\{n_{0},n_{1},....$\\
$,n_{k}\}.$ If index $[B]$ of $B$ with respect to $V$ is equal
to the multi index of $V$, then $B$ becomes a M-basis. $\newline$\\
$\mathbf{Proof.}$ Let us assume that $\omega\geq n_{0}>n_{1}>...>n_{k}\geq0.$
Then $\{\theta\}\subsetneqq V_{n_{0}}\subsetneqq V_{n_{1}}\subsetneqq V_{n_{2}}\subsetneqq....\subsetneqq V_{n_{k}}=X$.
Suppose that $[B]_{M}=\{r_{i}:i=0,1,2,...k\}.$ Then $dim\; V_{n_{i}}=\stackrel[j=0]{i}{\sum}r_{j}=\mid B\cap V_{n_{i}}\mid,$
for all $i=0,1,2,...,k.$ Hence, $B\cap V_{n_{i}}$ becomes a basis
for $V_{n_{i}}$ for each $i=0,1,2,..,k.$ Thus by $Proposition\;5.5,$
$B$ is a M-basis for $V.$$\newline$\\

\section{Dimension of multi vector space}

$\mathbf{Definition\:6.1}$ We define the dimension of a multi vector
space $V$ over $X$ by $dim(V)=\underset{B\: a\: base\: for\: X}{sup}\left(\underset{x\in B}{\sum}C_{V}(x)\right)$.\\
Clearly $dim$ is a function from the set of all multi vector spaces
to $\mathbb{N}.$ $\newline$\\
$\mathbf{Proposition\:6.2}$ Let $V\in MV(X)$ where $dim\: X=m<\infty.$
Then if $B$ is a M-basis for $V$ and $B^{*}$ is any basis for $X$
then \\
$\underset{x\in B^{*}}{\sum}C_{V}(x)\leq\underset{x\in B}{\sum}C_{V}(x).$
$\newline$\\
$\mathbf{Proposition\:6.3}$ If $V$ is a multi vector space over
a finite dimensional vector space $X$, then $dim(V)=\underset{x\in B}{\sum}C_{V}(x),$
where $B$ is any M-basis for $V.$ $\newline$\\
$\mathbf{Note\;6.4}$ If $V$ is a multi vector space over a finite
dimensional vector space $X$, then $dim(V)$ is independent of M-basis
for $V.$ It follows from $Proposition\;5.5$ and $Proposition\;5.7.$
$\newline$\\
$\mathbf{Proposition\;6.5}$ Let $X$ be any finite dimensional vector
space and $V,W\in MV(X)$ such that $C_{V}(\theta)\geq sup[C_{W}(X\setminus\{\theta\})]$
and $C_{W}(\theta)\geq sup[C_{V}(X\setminus\{\theta\})]$. Then there
exists a basis $B$ for $X$ which is also a M-basis for $V,$ $W$,
$V\cap W$ and $V+W.$ In addition, if $A_{1}=\{x\in X:C_{V}(x)<C_{W}(x)\},$
$A_{2}=X\setminus A_{1},$ then for all $v\in B\cap A_{1},$\\
$(C_{V\cap W})(v)=C_{V}(v)$ and $C_{V+W}(v)=C_{W}(v)$\\
and for all $v\in B\cap A_{2},$\\
$(C_{V\cap W})(v)=C_{W}(v)$ and $C_{V+W}(v)=C_{V}(v).$$\newline$\\
$\mathbf{Proof.}$ We prove this by induction on $dim\; X.$ In case
$dim\; X=1$ the statement is clearly true.\\
Now suppose that the theorem is true for all the multi vector space
with dimension of the underlying vector space equal to n.\\
Let $V$ and $W$ be two multi vector spaces over $X$ with $dim\; X=n+1>1.$
Let $B_{1}=\{v_{i}\}_{i=1}^{n+1}$ be any M-basis for $V.$ We may
assume that $C_{V}(v_{1})\leq C_{V}(v_{i})$ for all $i=\{2,3,...,n+1\}.$
Let $H=\prec\{v_{i}\}_{i=2}^{n+1}\succ.$ Since $n+1>1,$ $H\neq\{\theta\}.$
Clearly $dim\; H=n.$ Define the following two multi vector spaces:
$V_{1}$ with count function $C_{V_{1}}=C_{V}\mid_{H}$ and $W_{1}$
with the count function $C_{W_{1}}=C_{W}\mid_{H}$. By inductive hypothesis
the exists a basis $B^{*}$ for $H$ which is also a M-basis for $V_{1},$
$W_{1},$ $V_{1}\cap W_{1}$ and $V_{1}+W_{1}.$ Also for all $v\in B^{*}\cap A_{1},$\\
$(C_{V_{1}\cap W_{1}})(v)=C_{V_{1}}(v)$ and $C_{V_{1}+W_{1}}(v)=C_{W_{1}}(v)$\\
and for all $v\in B^{*}\cap A_{2},$\\
$(C_{V_{1}\cap W_{1}})(v)=C_{W_{1}}(v)$ and $C_{V_{1}+W_{1}}(v)=C_{V_{1}}(v).$
\\
We shall now show that $B^{*}$ can be extended to $B$ such that
$B$ is a M-basis for $V,$ $W$, $V\cap W$ and $V+W.$ Furthermore,
for all $v\in B\cap A_{1},$\\
$(C_{V\cap W})(v)=C_{V}(v)$ and $C_{V+W}(v)=C_{W}(v)$\\
and for all $v\in B\cap A_{2},$\\
$(C_{V\cap W})(v)=C_{W}(v)$ and $C_{V+W}(v)=C_{V}(v).$\\
\\
\textbf{Step - 1:} First we have to show that for all $x\in H,$\\
$C_{(V+W)}\mid_{H}(x)=C_{V_{1}+W_{1}}(x).$..........(1)\\
Let $x\in H.$ Then we have\\
$C_{(V+W)}\mid_{H}(x)=sup\{C_{V}(x_{1})\wedge C_{W}(x-x_{1}):\; x_{1}\in X\}$\\
$=sup\{C_{V}(x_{1})\wedge C_{W}(x-x_{1}):x_{1}\in H\}\vee sup\{C_{V}(x_{2})\wedge C_{W}(x-x_{2}):\; x_{2}\in X\setminus H\}.$........(2)\\
If $x\in H\setminus\{\theta\},$ we have\\
$C_{V}(x)\wedge C_{W}(x-x)=C_{V}(x)\wedge C_{W}(\theta)\leq sup\{C_{V}(x_{1})\wedge C_{W}(x-x_{1}):\; x_{1}\in H\},$
\\
$C_{V}(\theta)\wedge C_{W}(x-\theta)=C_{V}(\theta)\wedge C_{W}(x)\leq sup\{C_{V}(x_{1})\wedge C_{W}(x-x_{1}):\; x_{1}\in H\}.$
\\
Since $C_{V}(\theta)\geq sup[C_{W}(H\setminus\{\theta\})]$ and $C_{W}(\theta)\geq sup[C_{V}(H\setminus\{\theta\})],$\\
$C_{V}(x)\wedge C_{W}(\theta)=C_{V}(x)$ and $C_{V}(\theta)\wedge C_{W}(x)=C_{W}(x),$\\
and this leads to the following inequality:\\
$C_{V}(x)\vee C_{W}(x)\leq sup\{C_{V}(x_{1})\wedge C_{W}(x-x_{1}):\; x_{1}\in H\}.$.........(3)\\
Suppose that $sup\{C_{V}(x_{1})\wedge C_{W}(x-x_{1}):\; x_{1}\in H\}<sup\{C_{V}(x_{2})\wedge C_{W}(x-x_{2}):\; x_{2}\in X\setminus H\}.$...........(4)\\
This means that there exists $x^{\prime}\in X\setminus H$ such that
\\
$sup\{C_{V}(x_{1})\wedge C_{W}(x-x_{1}):\; x_{1}\in H\}<C_{V}(x^{\prime})\wedge C_{W}(x-x^{\prime}).$\\
In view of (2) we must have\\
$C_{V}(x)\vee C_{W}(x)<C_{V}(x^{\prime})\wedge C_{W}(x-x^{\prime}).$.......(5)\\
Since $x^{\prime}\in X\setminus H$ and $C_{V}(X\setminus H)=C_{V}(v_{1})\leq C_{V}(v_{i})$
for all $i\in\{2,3,...,n+1\},$ we must have $C_{V}(x)\geq C_{V}(x^{\prime}).$
Thus $(4)$ becomes $C_{V}(x)\vee C_{W}(x)<C_{V}(x)\wedge C_{W}(x-x^{\prime}).$
It is not possible (Using the properties of $\vee,\wedge$ and $<$).
This means that our assumption (3) is false. Therefore we must have
\\
$sup\{C_{V}(x_{1})\wedge C_{W}(x-x_{1}):\; x_{1}\in H\}\geq sup\{C_{V}(x_{2})\wedge C_{W}(x-x_{2}):\; x_{2}\in X\setminus H\}......(6)$\\
 If $x=\theta,$ $C_{(V+W)}\mid_{H}(\theta)=C_{V}(\theta)\wedge C_{W}(\theta)=sup\{C_{V}(x_{1})\wedge C_{W}(\theta-x_{1}):x_{1}\in H\}\geq sup\{C_{V}(x_{2})\wedge C_{W}(x-x_{2}):\; x_{2}\in X\setminus H\}....(7)$.
\\
Using equations $(6)$ and $(7)$,we have for all $x\in H,$\\
 $sup\{C_{V}(x_{1})\wedge C_{W}(x-x_{1}):\; x_{1}\in H\}\geq sup\{C_{V}(x_{2})\wedge C_{W}(x-x_{2}):\; x_{2}\in X\setminus H\},$\\
From $(1),$ we have for all $x\in H,$ $C_{(V+W)}\mid_{H}(x)=sup\{C_{V}(x_{1})\wedge C_{W}(x-x_{1}):\; x_{1}\in H\}$\\
$=sup\{C_{V}\mid_{H}(x_{1})\wedge C_{W}\mid_{H}(x-x_{1}):\; x_{1}\in H\}$\\
$=C_{V_{1}+W_{1}}(x).$\\
This establishes (1).\\
\\
Since $B^{*}$ is a M-basis of $V_{1}+W_{1},$$(1)$ implies that
$B^{*}$ is multi linearly independent in $V+W.$\\
\\
\textbf{Step - 2:} Let $v^{*}\in X\setminus H$ such that $C_{W}(v^{*})=sup[C_{W}(X\setminus H)].$
By $Lemma\;5.2$ and $Lemma\;5.3$, $v^{*}$ is an extension of M-basis
$B^{*}$ for $W_{1}$ to $B=B^{*}\cup\{v^{*}\}$ a M-basis for $W.$\\
\\
\textbf{Step - 3:} Since $C_{V}(X\setminus H)=C_{V}(v_{1}),$ $C_{V}(v_{1})=C_{V}(v^{*})$
and then $v^{*}$is also an extension of M-basis $B^{*}$ for $V_{1}$
to $B$ a M-basis for $V.$\\
\\
\textbf{Step - 4:} Now we shall show that $v^{*}$ is an extension
of M-basis $B^{*}$ for $V_{1}\cap W_{1}$ to $B$ a M-basis for $V\cap W.$
If $v^{*}\in A_{1}$, then $(C_{V}\wedge C_{W})(A_{1}\cap(X\setminus H))=C_{V}(v^{*})$
and for all $z\in A_{2}\cap(X\setminus H),$ $(C_{V}\wedge C_{W})(z)\leq C_{V}(v^{*})$,
by definition of $A_{2}.$ \\
From this we may conclude that if $v^{*}\in A_{1}$ then \\
$(C_{V}\wedge C_{W})(v^{*})=sup[(C_{V}\wedge C_{W})(X\setminus H)].$\\
If $v^{*}\in A_{2}$ then $C_{W}(v^{*})\leq C_{V}(v^{*}).$ Since
$C_{W}(v^{*})=sup[C_{W}(X\setminus H)]$ and $C_{V}$ is constant
on $X\setminus H$ we must have $A_{1}\cap(X\setminus H)=\phi.$ Therefore
we have that if $v^{*}\in A_{2}$, then $(C_{V}\wedge C_{W})(v^{*})=sup[(C_{V}\wedge C_{W})(X\setminus H)].$
By $Lemma\;5.3,$ we may now conclude that $v^{*}$ extends M-basis
$B^{*}$ for $V_{1}\cap W_{1}$ to $B$ a M-basis for $V\cap W.$\\
\\
\textbf{Step - 5:} Now we shall show that $v^{*}$ is also an extension
of $B^{*}$ a M-basis for $V_{1}+W_{1}$ to $B$ a M-basis for $V+W.$
Suppose that there exists $z\in X\setminus H$ such that $C_{(V+W)}(v^{*})<C_{(V+W)}(z).$
Clearly vector $z$ can be written in the form $z=a(v^{*}+v)$ where
$a\neq0$ and $v\in H.$ Therefore we have\\
$C_{(V+W)}(v^{*})<C_{(V+W)}(z)=C_{(V+W)}(a(v^{*}+v))=C_{(V+W)}(v^{*}+v).$
\\
This means that there exists $x_{1}\in X$ such that for all $x^{\prime}\in X,$\\
$C_{V}(x^{\prime})\wedge C_{W}(v^{*}-x^{\prime})<C_{V}(x_{1})\wedge C_{W}(v^{*}+v-x_{1}).$
............(8)\\
In particular this is true for $x^{\prime}=\theta,$ i.e. \\
$C_{V}(\theta)\wedge C_{W}(v^{*})<C_{V}(x_{1})\wedge C_{W}(v^{*}+v-x_{1}).$\\
But since $C_{V}(\theta)\geq sup[C_{W}(X\setminus\{\theta\})]$ we
have \\
$C_{W}(v^{*})<C_{V}(x_{1})\wedge C_{W}(v^{*}+v-x_{1}).$........(9)\\
If $x_{1}\in H$ then since $v\in H$ we must have $v-x_{1}\in H.$
Again $v^{*}\in X\setminus H.$ So, by $Lemma\;5.2,$ $C_{W}(v^{*}+v-x_{1})=C_{W}(v^{*})\wedge C_{W}(v-x_{1})$
and so (9) becomes $C_{W}(v^{*})<C_{V}(x_{1})\wedge C_{W}(v^{*})\wedge C_{W}(v-x_{1}),$
which is impossible. Thus $x_{1}\in X\setminus H.$ Let $x^{\prime}=v^{*}$
in (5). Since $C_{W}(\theta)\geq sup[C_{V}(X\setminus\{\theta\})]$
we have\\
$C_{V}(v^{*})<C_{V}(x_{1})\wedge C_{W}(v^{*}+v-x_{1}).....(10)$\\
Recall that $C_{V}(X\setminus H)=C_{V}(v_{1})$ and thus $C_{V}(v_{1})=C_{V}(v^{*})=C_{V}(x_{1})$,
as $x_{1}\in X\setminus H.$ This again means that the inequality
(10) is false. This means that for all $z\in X\setminus H,$ $C_{(V+W)}(v^{*})\geq C_{(V+W)}(z)$.
Therefore by $Lemma\;5.3,$ $v^{*}$ is an extension of $B^{*}$ a
M-basis for $V_{1}+W_{1}$ to $B$ a M-basis for $V+W.$\\
\\
\textbf{Step - 6:} Now we shall show that if $v^{*}\in A_{1}$ then
$C_{V+W}(v^{*})=C_{W}(v^{*})$ and if $v^{*}\in A_{2}$ then $C_{V+W}(v^{*})=C_{V}(v^{*}).$
From the definition we have:\\
$C_{(V+W)}(v^{*})=sup\{C_{V}(x_{1})\wedge C_{W}(v^{*}-x_{1}):\; x_{1}\in X\}.$\\
Let $x^{\prime}$ be such that \\
$sup\{C_{V}(x_{1})\wedge C_{W}(v^{*}-x_{1}):\; x_{1}\in X\}=C_{V}(x^{\prime})\wedge C_{W}(v^{*}-x^{\prime}).$\\
By substituting $x_{1}=\theta$ and then $x_{1}=v^{*}$ and recalling
that $C_{V}(\theta)\geq sup[C_{W}(X\setminus\{\theta\})]$ and $C_{W}(\theta)\geq sup[C_{V}(X\setminus\{\theta\})]$,
we obtain\\
$C_{V}(v^{*})\vee C_{W}(v^{*})\leq C_{V}(x^{\prime})\wedge C_{W}(v^{*}-x^{\prime}).$\\
Suppose that \\
$C_{V}(v^{*})\vee C_{W}(v^{*})<C_{V}(x^{\prime})\wedge C_{W}(v^{*}-x^{\prime}).$
.........(11)\\
If $x^{\prime}\in H$ then by $Lemma\;5.2$ (as $B=B^{*}\cup\{v^{*}\}$
is a M-basis for $W$), (11) becomes\\
$C_{V}(v^{*})\vee C_{W}(v^{*})<C_{V}(x^{\prime})\wedge C_{W}(v^{*})\wedge C_{W}(x^{\prime}).$
\\
This is never true, and thus $x^{\prime}\in X\setminus H.$ But now
since $C_{V}(v^{*})=C_{V}(x^{\prime})$ the inequality (11) never
holds, and so,\\
$C_{V}(v^{*})\vee C_{W}(v^{*})=C_{V}(x^{\prime})\wedge C_{W}(v^{*}-x^{\prime})=C_{(V+W)}(v^{*}).$
........(12)\\
From equation (12) , we have $v^{*}\in A_{1}$ then $C_{V+W}(v^{*})=C_{W}(v^{*})$
and if $v^{*}\in A_{2}$ then $C_{V+W}(v^{*})=C_{V}(v^{*}).$ \\
This completes the proof.$\newline$ \\
$\mathbf{Corollary\;6.6}$ If $V$ and $W$ are two multi vector spaces
over $X$ such that the dimension of $X$ is finite and $C_{V}(\theta)\geq sup[C_{W}(X\setminus\{\theta\})]$
and $C_{W}(\theta)\geq sup[C_{V}(X\setminus\{\theta\})]$, then\\
$dim(V+W)=dim\; V+dim\; W-dim\;(V\cap W)$. $\newline$ \\
$\mathbf{Example\;6.7}$ Suppose $X=\mathbb{R}^{2},\;\omega=6.$ Define
two multi vector spaces $V$ and $W$ with count functions $C_{V}$
and $C_{W}$ respectively as follows:\\
$C_{V}((0,0))=5$; $C_{V}((0,\mathbb{R}\setminus\{0\}))=3;$ $C_{V}(X\setminus\mathbb{R})=1$,\\
$C_{W}((0,0))=6;$ $C_{W}(\{(x,x):x\in\mathbb{R}\setminus\{0\}\})=2;$
$C_{W}(X\setminus\{(x,x):x\in\mathbb{R}\})=1.$\\
It is easily checked that $V$ and $W$ are multi vector spaces and
$C_{V}(\theta)\geq$\\
$sup[C_{W}(X\setminus\{\theta\})]$ and $C_{W}(\theta)\geq sup[C_{V}(X\setminus\{\theta\})]$.
It is also easy to check that \\
$C_{V\cap W}((0,0))=5$, $C_{V\cap W}(\{(x,x):x\in\mathbb{R}\setminus\{0\}\})=1$,
$C_{V\cap W}(X\setminus\{(x,x):x\in\mathbb{R}\})=1,$ $C_{V+W}((0,0))=5$;
$C_{V+W}((0,\mathbb{R}\setminus\{0\}))=3;$ $C_{V+W}(X\setminus(0,\mathbb{R}))=2$
and $B=\{(0,1),(1,1)\}$ is a M-basis for $V,$ $W,$ $V\cap W$ and
$V+W.$ Thus\\
$dim\;(V+W)=3+2=5,$ $dim(V\cap W)=1+1=2,$\\
$dim\; V=3+1=4$, $dim\; W=2+1=3$.\\
So, $dim\; V+dim\; W-dim\;(V\cap W)=4+3-2=5=dim\;(V+W).$ $\newline$
\\
$\mathbf{Definition\:6.8}$ Let $V$ be a multi vector space over
$X$ and $f:X\rightarrow Y$ be a linear map. Then we define $f(V)$
as\\
$C_{f(V)}(x)=\begin{cases}
sup\{C_{V}(z):z\in f^{-1}(x)\} & \: if\: f^{-1}(x)\neq\phi\\
0 & otherwise
\end{cases}$\\
and $\tilde{kerf}=(kerf,C_{V}\mid_{kerf}),\:\tilde{imf}=(imf,C_{V}\mid_{imf})$.$\newline$
\\
$\mathbf{Proposition\:6.9}$ If $V$ be a multi vector space over
$X$ where $dim\; X$ is finite and $f:X\rightarrow Y$ is a linear
map, then \\
$dim(\tilde{kerf})+dim(\tilde{imf})=dim(V).$ $\newline$ \\
$Proof.$ Suppose that $kerf\neq\{\theta\}.$ If $kerf=\{\theta\}$
then the proof is similar. Now let $B_{Kerf}$ be a M-basis for $\tilde{kerf}$
and $B_{Ex}$ be an extension of $B_{Ker}$ to a M-basis for $V$
(this is clearly possible by repeated application of $Lemma\:5.3$).
Then $B_{Kerf}\cup B_{Ex}=B$ is M-basis for $V$ and $B_{Kerf}\cap B_{Ex}=\phi.$
\\
We first show that $f(B_{Ex})=B_{Im}$ is a M-basis for $\tilde{imf}.$
Clearly $B_{Im}$ is a basis for $imf.$ Let $v_{1},v_{2},..,v_{k}\in B_{Ex}$
and $a_{1},...,a_{k}\in\mathbb{R}$ not all zero. By definition we
have \\
$C_{f(V)}(\stackrel[i=1]{k}{\sum}a_{i}f(v_{i}))$\\
$=\begin{cases}
sup\{C_{V}(x):x\in f^{-1}(\stackrel[i=1]{k}{\sum}a_{i}f(v_{i}))\} & if\: f^{-1}\left(\stackrel[i=1]{k}{\sum}a_{i}f(v_{i})\right)\neq\phi\\
0 & otherwise
\end{cases}$\\
Since $\stackrel[i=1]{k}{\sum}a_{i}f(v_{i})\in imf$ we have \\
$C_{f(V)}(\stackrel[i=1]{k}{\sum}a_{i}f(v_{i}))=sup\{C_{V}(x):x\in f^{-1}(\stackrel[i=1]{k}{\sum}a_{i}f(v_{i}))\}.$\\
By linearity of $f$ and by the property of $f^{-1}$ we get \\
$C_{f(V)}(\stackrel[i=1]{k}{\sum}a_{i}f(v_{i}))=sup\{C_{V}(x):x\in kerf+\stackrel[i=1]{k}{\sum}a_{i}v_{i}\}.$\\
If $x\in kerf$ then $x=\theta$ or $x=\stackrel[i=1]{p}{\sum}b_{i}u_{i},u_{i}\in B_{kerf}$
where not all $b_{i}$ are zero; so if $x\in kerf+\stackrel[i=1]{k}{\sum}a_{i}v_{i}$
then either $C_{V}(x)=C_{V}(\theta+\stackrel[i=1]{k}{\sum}a_{i}v_{i})$
or $C_{V}(x)=C_{V}(\stackrel[i=1]{p}{\sum}b_{i}u_{i}+\stackrel[i=1]{k}{\sum}a_{i}v_{i})$
and thus\\
$C_{V}(x)=min\left(\stackrel[i=1]{p}{\wedge}C_{V}(b_{i}u_{i}),\stackrel[i=1]{k}{\wedge}C_{V}(a_{i}v_{i})\right),$
{[} As $u_{i}$and $v_{i}$ are M-basis element of $V${]}\\
which is clearly smaller than or equal to $C_{V}(\stackrel[i=1]{k}{\sum}a_{i}v_{i}).$
Thus\\
$C_{f(V)}(\stackrel[i=1]{k}{\sum}a_{i}f(v_{i}))=sup\{C_{V}(x):x\in kerf+\stackrel[i=1]{k}{\sum}a_{i}v_{i}\}=C_{V}(\stackrel[i=1]{k}{\sum}a_{i}v_{i})=\stackrel[i=1]{k}{\wedge}C_{V}(a_{i}v_{i}).$\\
By the same argument we get that $C_{f(V)}(f(v_{i}))=C_{V}(v_{i}).$
Thus $C_{f(V)}(\stackrel[i=1]{k}{\sum}a_{i}f(v_{i}))=\stackrel[i=1]{k}{\wedge}C_{f(V)}(a_{i}v_{i}).$\\
and therefore $B_{Im}$ is a M-basis for $\tilde{imf}.$\\
Now by definition of multi dimension we get \\
$dim(V)=\underset{v\in B_{Ker}\cup B_{Ex}}{\sum}C_{V}(v)=\underset{v\in B_{Ker}}{\sum}C_{V}(v)+\underset{v\in B_{Ex}}{\sum}C_{V}(v)$.\\
But by the above we have if $z\in\text{\ensuremath{\prec}}B_{Ex}\text{\ensuremath{\succ}}$,
then $C_{f(V)}(f(z))=C_{V}(z),$ and thus \\
$dim(V)=\underset{v\in B_{Ker}}{\sum}C_{V}(v)+\underset{v\in B_{Ex}}{\sum}C_{f(V)}(f(v))$\\
$=\underset{v\in B_{Ker}}{\sum}C_{V}(v)+\underset{v\in B_{Im}}{\sum}C_{f(V)}(v)$\\
$=dim(\tilde{kerf})+dim(\tilde{imf}).$

\section{Conclusion}

There is a future scope of study of infinite dimensional multi vector
space and behavior of linear operators in multi vector space context.


\begin{thebibliography}{10}
\bibitem{bli} W. D. Blizard, Multiset theory, Notre Dame J. Formal
Logic, 30 (1989) 36-66. 

\bibitem{chakraborty} K. Chakraborty, On bags and fuzzy bags, Adv.
Soft Comput. Techniq. Appl. 25 (2000) 201-21

\bibitem{chiney} M. Chiney, S. K. Samanta, Multi vector space, Annals
of Fuzzy Mathematics and Informatics 13 (5) 553-562.

\bibitem{sujoy} S. Das, P. Majumdar and S. K. Samanta, On soft linear
spaces and soft normed linear spaces, Annals of Fuzzy Mathematics
and Informatics 9 (1) (2015) 91-109.

\bibitem{delago} M. Delgado, M. J. Martin Bautista, D. Sanchez and
M. A. Vila, An extended characterization of fuzzy bags, Int. J. Intell.
Syst. 24 (2009) 706-721.

\bibitem{delago2} M. Delgado, M. D. Ruiz and D. Sanchez, Pattern
extraction from bag data bases, Internat. J. Uncertain. Fuzziness
Knowledge-Based Systems 16 (2008) 475-494.

\bibitem{girish} K. P. Girish, S. J. John, Relations and functions
in multiset context, Inform. Sci. 179(6)(2009) 758-768.

\bibitem{girish-1} K. P. Girish, S. J. John, General relations between
partially ordered multisets and their chains and antichains, Math.
Commun. 14(2) (2009) 193-206.

\bibitem{hickman} J. L. Hickman, A note on the concept of multiset,
Bull. Austral. Math. Soc. 22(2) (1980) 211-217.

\bibitem{katsaras} A. K. Katsaras, D. B. Liu, Fuzzy vector spaces
and fuzzy topological vector spaces, J. Math. Anal. Appl. 58(1977)
135-146.

\bibitem{klausner} A. Klausner, N. Goodman, Multirelations-Semantics
and languages, in: Proceedings of the 11th Conference on Very Large
Data Bases VLDB'85 (1985) 251-258.

\bibitem{kosters} W. A. Kosters, J. F. Laros, Metrics for mining
multisets, in: 27th SGAI International Conference on Innovative Techniques
and Applications of Artificial Intelligence (2007) 293-303.

\bibitem{li} B. Li, W. Peizhang and L. Xihui, Fuzzy bags with set-valued
statistics, Comput. Math. Appl. 15 (1988) 3-39. 

\bibitem{lubzonok} P. Lubzonok, Fuzzy vector spaces, Fuzzy Sets and
Systems 38 (1990) 329-343.

\bibitem{majumdar} P. Majumdar and S. K. Samanta, Soft multisets,
J. Math. Comput. Sci. 2(6) (2012) 1700-1711.

\bibitem{miyamoto} S. Miyamoto, Operations for real-valued bags and
bag relations, in ISEA-EUSFLAT, (2009) 612-617.

\bibitem{miyamoto1} S. Miyamoto, Two generalizations of multisets,
in: M. Inuiguchi, S. Tsumito eds., Rough Set theory \& Granular computing,
Springer, (2003) 59-68.

\bibitem{mumick} I. S. Mumick, H. Pirahesh and R. Ramakrishnan, The
magic of duplicates and aggregates, in: Proceedings of the 16th Conference
on Very Large Data Bases VLDB'90 (1990) 264-277.

\bibitem{nazmul} Sk. Nazmul, S. K. Samanta, Multisets and Multigroups,
Annals of fuzzy mathematics and informatics 6(3) (2013) 643-656.

\bibitem{paun} G. Paun and M. J. Perez-Jimenez, Membrane computing:
brief introduction, recent results and applications, Bio Systems 85
(2006) 11-22.

\bibitem{pradhan} R. Pradhan , M. Pal, Intuitionistic Fuzzy Linear
Transformations, Annals of Pure and Applied Mathematics Vol. 1 (1)
(2012) 57-68.

\bibitem{shi} F. G. Shi, C E Huang, Fuzzy bases and the fuzzy dimension
of fuzzy vector spaces, Math. Commun., 15 (2) (2010) 303-310.

\bibitem{yager} R. R. Yager, On the theory of bags, Internat. J.
Gen. Systems 13(1) (1987) 23-37.\end{thebibliography}
\end{document}